 %&AMS-TeX
%\input amstex.tex
\documentstyle{amsppt}
\magnification=1200
\hsize=150truemm
\vsize=224.4truemm
\hoffset=4.8truemm
\voffset=12truemm
 
%\NoBlackBoxes
\NoRunningHeads

\define\1{{\bold 1}}%                \1          1 gras
\define\C{{\bold C}}%                \C          C gras
 
%                \Q          Q gras
%                \R          R gras
\define\Z{{\bold Z}}%                \Z          Z gras

\let\thm\proclaim%             \thm        \proclaim
\let\fthm\endproclaim%         \fthm       \endproclaim    
 
%%%%%%%%    numerotation des sections et enonces %%%%%%%%%%%%
\newcount\tagno
\newcount\secno
\newcount\subsecno
\newcount\stno
\global\subsecno=1
\global\tagno=0
\define\ntag{\global\advance\tagno by 1\tag{\the\tagno}}

\define\sta{\ %\uppercase\expandafter{\romannumeral
{\the\secno}.\the\stno
\global\advance\stno by 1}

\define\stas{\the\stno
\global\advance\stno by 1}

\define\sect{\global\advance\secno by 1
\global\subsecno=1\global\stno=1\
%\uppercase\expandafter{\romannumeral\the\secno}. }
{\the\secno}. }

\def\nom#1{\edef#1{{\the\secno}.\the\stno}}
\def\inom#1{\edef#1{\the\stno}}
\def\eqnom#1{\edef#1{(\the\tagno)}}

%%%%%%%%%%%%%%  macros pour la bibliographie  %%%%%%%%%%%%%%%%%%%%%%
\newcount\refno
\global\refno=0

\def\nextref#1{
      \global\advance\refno by 1
      \xdef#1{\the\refno}}

\def\bref {\ref\global\advance\refno by 1\key{\the\refno}}

%%%%%%%%% bibliographie  %%%%%%%%%%%
\nextref\AH
\nextref\BD
\nextref\BL
\nextref\BO
\nextref\BR
\nextref\CZ
\nextref\D
\nextref\DU
\nextref\GR
\nextref\GH
\nextref\HU
\nextref\KL
\nextref\LE
\nextref\MA
\nextref\MC
\nextref\PA
\nextref\RO
\nextref\SZ
\nextref\WI

%%%%%%%%% en tete%%%%%%%
\topmatter

 %\abstract
 %\endabstract
\title 
Around brody lemma
 \endtitle

\author  Julien Duval\footnote""{Laboratoire de Math\'ematiques, Universit\'e Paris-Sud, 91405 Orsay cedex, France \newline julien.duval\@math.u-psud.fr\newline}
\footnote""{Keywords : Brody curve, Ahlfors current, hyperbolicity \newline AMS class. : 32Q45\newline}
\endauthor
 
\endtopmatter 

\document

\subhead 1. Brody lemma \endsubhead

\null
Let $X$ be a compact complex manifold. An {\it entire curve} in $X$ is a non constant holomorphic map $f: \C \to X$. It is a {\it Brody curve} if its derivative $\Vert f'\Vert$ is bounded, where the norm is computed with respect to the standard metric on $\C$ and a given metric on $X$. Brody curves arise naturally as limits of sequences of larger and larger holomorphic discs, thanks to Brody lemma [\BR].
\thm {Brody lemma} Let $f_n : D\to X$ a sequence of holomorphic maps from the unit disc to a compact complex manifold. Suppose that $\Vert f_n'(0)\Vert$ blows up. Then there exist reparametrizations $r_n$ such that $f_n \circ r_n$ converges toward a Brody curve, after extracting a subsequence. \fthm

\demo {Proof} We may suppose $f_n$ smooth up to the boundary. Denote by $\delta(z)$ the distance from $z$ to $\partial D$. As the function $\delta \Vert f_n'\Vert$ vanishes on $\partial D$ it reaches its maximum inside $D$, say at $a_n$. This is where we will reparametrize $f_n$. Let $d_n$ be the disc $D(a_n, \frac {\delta(a_n)}{2})$. We have $ \Vert\delta f'_n\Vert _{d_n}\leq  \delta(a_n)\Vert f_n'(a_n)\Vert $. But $\delta\geq \frac {\delta(a_n)}{2}$ on $d_n$ so
$\Vert f'_n\Vert _{d_n}\leq 2\Vert f_n'(a_n)\Vert $. 

Define the reparametrization by $r_n(z)=a_n+ \frac {z} {\Vert f_n'(a_n)\Vert}$. let $D_n$ be the preimage of $d_n$ by $r_n$. Its radius blows up as $\delta(a_n)\Vert f_n'(a_n)\Vert \geq \Vert f_n'(0)\Vert$. We may suppose that $D_n$ increases toward $\C$ after extracting. Moreover $\Vert (f_n \circ r_n)' \Vert_{D_n} \leq 2$ and $\Vert (f_n\circ r_n)'(0)\Vert =1$. By Ascoli theorem we may extract a subsequence of $f_n\circ r_n$ which converges toward a holomorphic map $f : \C \to X$ such that $\Vert f'\Vert_{\C}\leq 2$ and $\Vert f'(0)\Vert =1$. It is a Brody curve.
\enddemo

As a consequence we get a characterization of Kobayashi-hyperbolicity. 

\null\noindent
Recall that the {\it Kobayashi pseudometric} of $X$ at $p$ in the direction $v$ is $K(p,v)=\inf \{r>0 /\exists f:D\to X \text{ holomorphic, } f(0)=p,f'(0)=\frac v r\}$. It measures the size of the holomorphic discs passing through a point in a given direction. The manifold $X$ is {\it Kobayashi-hyperbolic} if its pseudometric is non degenerate. 

\thm {Criterion} A compact complex manifold $X$ is Kobayashi-hyperbolic if and only if there is no entire curve in it.
\fthm

\noindent Indeed the vanishing of $K(p,v)$ gives rise to a sequence of holomorphic discs $f_n : D \to X$ such that $f_n(0)=p$ and $f_n'(0)$ blows up in the direction of $v$. By Brody lemma we get a Brody curve in $X$. Conversely the Kobayashi pseudometric has to vanish along any entire curve.

\null

In the sequel we will say that $U\subset X$ is {\it hyperbolic} if $U$ does not contain any entire curve. Hyperbolicity is invariant by etale covering. Indeed if $U\to V$ is such a covering, an entire curve in $U$ may be pushed down in an entire curve in $V$. Conversely an entire curve in $V$ can be lifted to an entire curve in $U$. Another consequence of Brody lemma is that hyperbolicity is an open property.

\thm {Openness} Let $X$ be a compact complex manifold and $F\subset X$ a closed subset. If $F$ is hyperbolic then so is a small neighbourhood of it. \fthm

\noindent Indeed if not we would get an entire curve in each $\epsilon$-neighbourhood $F_\epsilon$ of $F$. This would allow us to construct a sequence of holomorphic discs $f_n : D \to F_{\frac 1 n}$ such that $\Vert f_n'(0)\Vert$ blows up, giving at the limit a Brody curve in $F$.

\null
We discuss now some examples of hyperbolic surfaces. We start with the simplest hyperbolic complement.
\thm {Green theorem [\GR]} Let $L$ be a collection of five lines in general position in $P^2(\C)$. Then $P^2(\C)\setminus L$ is hyperbolic. \fthm

Here general position means there is no triple point in the configuration.
\demo {Proof} Embedd $P^2(\C)$ into $P^4(\C)$ by $z\mapsto [l_1(z):...:l_5(z)]$ using the equations of the lines. Call $P$ its image and $P^*=P\setminus H$ the complement of the collection $H$ of coordinate hyperplanes in $P^4(\C)$. We want to prove the hyperbolicity of $P^*$. Let $F_n$ be the self-map of $P^4(\C)$ given by $z\mapsto [z_1^n:...:z_5^n]$. It induces an etale covering from $P^4(\C)\setminus H$ to itself, so the hyperbolicity of $P^*$ and $F_n^{-1}(P^*)$ are equivalent. The point is that $F_n^{-1}(P)$ converges toward a polyhedron whose hyperbolicity is easily checked by Liouville theorem. This will conclude by openness. 

\null
Let us make it precise. By general position $P$ avoids the coordinate lines $(z_i=z_j=z_k=0)$. So $P$ is contained in $X_{\epsilon}=\cap_{\{i,j,k\}}(\text{max}(\vert z_i\vert, \vert z_j\vert, \vert z_k \vert )\geq \epsilon \Vert z\Vert)$ for some small $\epsilon>0$. Here $\Vert z \Vert =\text{max} \vert z_i \vert$. Now $F_n^{-1}(X_{\epsilon})=X_{\epsilon^{\frac 1 n}}$ decreases toward the polyhedron $X=X_1$ when $n$ goes to infinity. Note that on $X$ $\Vert z \Vert$ is reached on one out of three arbitrary components of $z$, meaning that $\Vert z \Vert$ is always reached on three components at least. So $X$ is alternatively seen as a finite union of faces $X_{i,j,k}=(\vert z_i\vert=\vert z_j\vert =\vert z_k\vert= \Vert z \Vert)$. 

Let us check its hyperbolicity. If $f:\C \to P^4(\C), z\mapsto [f_1(z):...:f_5(
z)]$ is a holomorphic map such that $f(\C) \subset X$, it has to spend some time in one of the faces, say $X_{1,2,3}$. This implies by analytic continuation that $\vert f_1\vert =\vert f_2\vert = \vert f_3\vert $ everywhere. But $\Vert f \Vert =\text{max}(\vert f_1\vert,\vert f_2\vert, \vert f_3\vert)$ by definition of $X$. Hence $\vert f_1\vert$ dominates the other components everywhere, meaning that $f$ is bounded in the chart $(z_1=1)$ thus constant by Liouville theorem. Therefore $X$ does not contain any entire map.

 \enddemo
 
This argument is dynamical in essence as $X$ is nothing but an intermediate Julia set for $F_2$ which is known to attract backward iterates of a generic plane. It stems from the proof of Picard theorem (the hyperbolicity of $P^1(\C)\setminus 3$ points) by A. Ros [\RO] and works in any dimension [\BD]. 

\null
We focus now on examples of hyperbolic surfaces of low degree in $P^3(\C)$. This fits into the framework of Kobayashi conjecture which predicts that a generic surface of degree $\geq 5$ in $P^3(\C)$ is hyperbolic. It holds true for degree $\geq 18$ [\PA] but few examples of hyperbolic surfaces of smaller degree are known, according to the motto "the lower the degree the harder the hyperbolicity". Here we adapt a deformation method due to M. Zaidenberg [\SZ] to produce examples of degree $6$ by reduction to Green theorem. This also works for higher degrees and higher dimensions [\HU]. Note that remains open the question of finding a single hyperbolic quintic in $P^3(\C)$.

\null
We will use Brody lemma in the following form.
\thm {Sequences of entire curves} Let $X$ be a compact complex manifold. Then any sequence of entire curves in $X$ can be made converging toward an entire curve after reparametrization and extraction.\fthm
\noindent 
Indeed let $(f_n)$ be the sequence of entire curves. By translating we may suppose that $f_n'(0)$ does not vanish and by dilating that actually $\Vert f_n'(0)\Vert$ blows up. It remains to apply Brody lemma.

\null
We will also invoke the following fact.
\thm {Stability of intersections} Let $X$ be a complex manifold and $H\subset X$ an analytic hypersurface. Suppose that a sequence $(f_n)$ of entire curves in $X$ converges toward an entire curve $f$. If $f(\C)$ is not contained in $H$ then $f(\C)\cap H \subset lim  f_n(\C) \cap H$. \fthm 
\noindent
Indeed let $z$ be a point in $f^{-1}(H)$ and $h$ a local equation of $H$ near $f(z)$. As $h\circ f$ does not vanish identically near $z$, $h\circ f_n$ has to have a zero in any small neighbourhood of $z$ for large $n$ by Rouch\'e theorem. 

\null
We construct now our example of hyperbolic sextic as a suitable small deformation of a union of six planes (see also [\CZ], [\D] for other examples).

\thm {A hyperbolic sextic} Let $(P_i=(p_i=0))$ be a collection of six planes in general position in $P^3(\C)$. Then we can find a sextic $S=(s=0)$ such that the surface $\Sigma_\epsilon=(\Pi p_i=\epsilon s)$ is hyperbolic for $\epsilon \neq 0$ sufficiently small. \fthm
\noindent Here general position means there is no quadruple point in the configuration. Moreover $S$ will be in general position with respect to the $P_i$, in the sense that it will avoid the triple points of the configuration. Note that by definition $\Sigma_\epsilon\cap P_i\subset S$.

\demo{Proof} The first step reduces the problem to the hyperbolicity of complements. It is the heart of Zaidenberg method and goes as follows. If $\Sigma_{\epsilon_n}$ is not hyperbolic for $\epsilon_n$ going to zero we have entire curves $f_n : \C \to \Sigma_{\epsilon_n}$. By Brody lemma we get at the limit an entire curve $f : \C \to \Sigma_0=\cup P_i$. It lands inside one of the planes. 

\null
We analyze now its position with respect to the rest of them. The crucial remark is the following. If $f(\C)$ is not contained in $P_i$ then $f(\C)\cap P_i\subset S$. Indeed by stability of intersections $f(\C)\cap P_i \subset \text {lim}\ f_n (\C)\cap P_i \subset \text{lim}\ \Sigma_{\epsilon_n} \cap P_i \subset S$. We infer that $f(\C)$ cannot land into a double line of the configuration of planes. If it were the case $f(\C)$ would have to avoid the 4 triple points on the line by the remark and the general position of $S$, contradicting Picard theorem. 

\null
We end up with $f(\C)$ contained in one plane and avoiding the others except at points of $S$, again by the remark. So $f(\C)$ is in a complement of the form $P_i\setminus (\cup_{j\neq i}P_j\setminus S)$. Hence we are finished if we are able to find a sextic $S$ such that all these complements are hyperbolic.

\null
The second step consists in constructing this sextic. Note that the situation is similar to Green theorem. We have plane complements of five lines on which a few points are deleted.
To create $S$ we proceed by deformation in order to remove these points on more and more double lines. Our starting point is the collection of complements $P_i\setminus (\cup_{j\neq i}P_j)$ which are hyperbolic by Green theorem. 

\null
We want to remove points on a double line, say $D=P_1\cap P_2$, keeping the hyperbolicity. For this consider a sextic $S_0=(s_0=0)$ in general position with respect to the $P_i$ and deform it toward the union of the remaining planes by taking $S_1= 
(p_3p_4p_5^2p_6^2=\epsilon_0 s_0)$ for a small $\epsilon_0\neq 0$. This pushes on $D$ the points of the sextic toward the triple points. So the complement $P_1\setminus(\cup_{i\neq1}P_i\setminus (D\cap S_1))$ is close to $P_1\setminus(\cup_{i\neq1}P_i)$, hence hyperbolic by a suitable openness argument. 

\null
This is a particular case of the following lemma which we apply inductively to conclude.

\enddemo
\thm {Lemma} Let $\Delta_k$ be a collection of $k$ double lines, $D=P_{i_{1}}\cap P_{i_{2}}$ an extra one and $\Delta_{k+1}=\Delta_k\cup D$. Assume $S_k=(s_k=0)$ already constructed such that the complements 
$P_i\setminus (\cup_{j\neq i} P_j\setminus (\Delta_k\cap S_k))$ are hyperbolic.
 Then so are $P_i\setminus (\cup_{j\neq i} P_j\setminus (\Delta_{k+1}\cap S_{k+1}))$ where $S_{k+1}= 
( p_{i_{3}}p_{i_{4}}p^2_{i_{5}}p_{i_{6}}^2=\epsilon_k s_k)$ for a small $\epsilon_k\neq0$. \fthm

\noindent Note that $S_{k+1}$ is still in general position with respect to the $P_i$ if $S_k$ was. Remark also that the situation does not change on $\Delta_k$. We have $\Delta_k\cap S_k=\Delta_k\cap S_{k+1}$.

\demo{Proof of the lemma} Take $D=P_1\cap P_2$ for simplicity. If we cannot find such an $\epsilon_k$, we have a sequence of sextics $S_{k+1,n}$ converging toward $P_3\cup P_4\cup P_5 \cup P_6$ and entire curves $f_n(\C)$ sitting in one of the corresponding complements, say $P_1\setminus (\cup_{j\geq 2} P_j\setminus ((\Delta_k\cap S_k)\cup (D\cap S_{k+1,n})))$. We get at the limit an entire curve $f(\C)$ in $P_1$. As before it cannot degenerate inside a double line. By stability of intersections, for $j\geq2$ we have $f(\C)\cap P_j \subset \text{lim}f_n(\C)\cap P_j\subset (\Delta_k\cap S_k)\cup \text{lim}\ D \cap S_{k+1,n}$. If $j\geq 3$ we infer that $f(\C) \cap P_j \subset \Delta_k \cap S_k$ as $P_j\cap D\cap S_{k+1,n}$ is empty by general position. Note now that $\text{lim}\ D\cap S_{k+1,n}$ consists in the triple points of $D$ hence sits in $\cup _{j\geq 3} P_j$. Then thanks to the previous case we also have $f\C) \cap P_2 \subset \Delta_k \cap S_k$. Therefore $f(\C)$ lands in $P_1\setminus (\cup_{j\geq 2} P_j\setminus (\Delta_k\cap S_k))$ contradicting the hypothesis.

\enddemo

\null

\subhead 2. A variant \endsubhead

\null
A drawback of Brody lemma is the lack of information about the location of the entire curve it produces. It might land far away from the points where the discs blow up. Here is a simple example due to J. Winkelmann (see also [\WI], and [\GH] for background on blow-ups).

\thm {Example} Let $A=\C^2/(\Z\oplus i\Z)^2$ be the standard torus and $\pi :\tilde A \to A$ the blow-up of $A$ at a point $p$. Take a dense injective line $L$ in $A$, say $L=(z_2=\lambda z_1)$ for $\lambda$ irrational. Consider the sequence of discs $f_n(D)$ on $L$ given by $f_n(z)=(nz, \lambda nz)$. Then the Brody curve $\tilde f(\C)$ produced by their strict transforms $\tilde f_n$ lands in the exceptional divisor $E$. \fthm

\noindent Indeed if $\tilde f(\C)$ is not contained in $E$ it projects down to a (non constant) Brody curve $f(\C)$ in $A$. It is linear by Liouville theorem (the lift of $f'$ is bounded in $\C^2$), and parallel to $L$ by construction. So it is again a dense injective line $L'$. In particular the derivative of $\pi^{-1}\vert_{L'}=\tilde f\circ f^{-1}$ is bounded. On the other hand consider an open cone of vertex $p$ transversal to the direction of $L'$. By density it cuts out a sequence of smaller and smaller discs $d_n$ on $L'$ converging toward $p$. But $\pi^{-1}(d_n)$ converges toward a non constant disc in $E$ (the basis of the cone) by definition of a blow-up. This is the contradiction.

\null

We want to address this problem by constructing entire curves where the discs accumulate in terms of area. To formulate this we introduce the notion of Ahlfors current. 

\null
We start by briefly recalling what a current is (see $[\GH]$ for more background on currents). Let $X$ be a compact complex manifold of dimension $n$ endowed with a hermitian metric. Denote by $\omega$ its area form. A {\it current} (of bidimension $(1,1)$) in $X$ is a continuous linear form on the space of differential forms (of bidegree $(1,1)$) of $X$. It is {\it positive} if it is non negative on positive forms, and {\it closed} if it vanishes on exact forms. Instances of such currents are currents of integration on complex curves or positive closed forms of bidegree $(n-1,n-1)$, given respectively by $\alpha \mapsto \int_C \alpha$ and $\alpha \mapsto \int_X \alpha\wedge \beta$. In general a positive current may be seen as a form of bidegree $(n-1,n-1)$ with measure coefficients. A useful fact is the compactness of the set of positive currents of mass $1$ (here the mass of $T$ is $<T,\omega>$ or $\int_X T\wedge \omega$ if we see $T$ as a form with measure coefficients). 

\null
Let now $f_n : D \to X$ be a sequence of holomorphic discs smooth up to the boundary. We make the standing assumption
$$\frac{l_n}{a_n}\to 0 \tag"(A)" $$
where $a_n=\text{area }(f_n(D))$ and $l_n=\text{length }(f_n(\partial D))$. In other words the boundaries of the discs become asymptotically negligible. By compactness we may suppose that the sequence of normalized currents of integration $\frac {[f_n(D)]}{a_n}$ converges toward a positive current $T$ of mass $1$. It is closed because of $(A)$. This is the {\it Ahlfors current} associated to the sequence of discs.

\null
Let us present our variant of Brody lemma.

\thm{Variant [\DU]} Let $T$ be an Ahlfors current in $X$ and $K\subset X$ a compact set charged by $T$. Then there exists an entire curve passing through $K$. \fthm

\noindent Here $K$ charged by $T$ means $\int_KT\wedge\omega>0$. In other words the proportion of area of $f_n(D)$ near $K$ remains bounded away from zero.
As in Brody lemma the entire curve $f:\C\to X$ is obtained by reparametrizing the holomorphic discs. By the very construction it cuts $K$ out on a set of positive area. In particular if $K$ is an analytic subset then $f(\C)\subset K$ by analytic continuation.

\null
As a consequence we get another characterization of hyperbolicity with a flavour of negative curvature (see also [\KL] for another proof).

\thm{Criterion} Let $X$ be a compact complex manifold. Then $X$ is hyperbolic if and only if its holomorphic discs satisfy a linear isoperimetric inequality.  \fthm

\noindent This means that there exists a constant $C$ such that for any holomorphic disc $f:D\to X$ smooth up to the boundary we have $$ \text{ area }(f(D))\leq C \text { length }(f(\partial D)).$$

\noindent Indeed if it does not hold we get a sequence of holomorphic discs satisfying $(A)$, therefore an Ahlfors current and so an entire curve in $X$ by the theorem. Conversely the presence of an entire curve forbids a linear isoperimetric inequality by Ahlfors lemma.

\thm{Ahlfors lemma [\AH]} Let $f: \C \to X$ be an entire curve. Then there exists $r_n\to+\infty$ such that the corresponding sequence of holomorphic discs $f_n :D\to X, z\mapsto f(r_nz)$ satisfies $(A)$.  \fthm
 
\demo{Proof} Let $D_r$ be the disc of radius $r$ centered at $0$ in $\C$. Denote by $l(r)$ the length of $f(\partial D_r)$ and $a(r)$ the area of $f(D_r)$. If $\lambda\vert dz\vert$ is the pull-back by $f$ of the metric on $X$ then in polar coordinates $$l(r)=\int_0^{2\pi}\lambda(r,\theta)\ r\ d\theta,\ \ \  a'(r)=\int_0^{2\pi}\lambda^2(r,\theta)\ r\ d\theta.$$ By Cauchy-Schwarz inequality we get $l^2(r)\leq 2\pi r \ a'(r)$. Dividing by $2\pi r\ a^2(r)$ and integrating we infer that $$\int_1^{+\infty}\frac {l^2(r)}{a^2(r)} \frac {dr}{2\pi r}\leq \int_1^{+\infty}\frac {a'(r)}{a^2(r)}dr \leq \frac 1 {a(1)} <+\infty.$$ This gives a sequence $r_n\to +\infty$ such that $\frac {l(r_n)}{a(r_n)}\to 0$.
\enddemo

Ahlfors lemma also gives a description of rational curves as entire curves of bounded area.

\thm{Rational curves} Let $f:\C\to X$ be an entire curve such that area($f(\C))<+\infty$. Then $f$ extends to a holomorphic map from $P^1(\C)$ to $X$, a rational curve.
\fthm

\noindent Indeed by Riemann removable singularity theorem it is enough to extend $f$ continuously at infinity. Because of the bounded area, applying Ahlfors lemma xe get a sequence $r_n$ such that length(f$(\partial D_{r_n}))\to 0$. We may suppose that $f(\partial D_{r_n})$ converges toward a point $p$. It suffices to show that the annuli $A_n=f(D_{r_n}\setminus\overline D_{r_{n-1}})$ also converge toward $p$. Note that area($A_n)\to 0$. If a point $q_n$ of $A_n$ remains far from $p$ say at distance $>\epsilon$, then for large $n$ $A_n\cap B(q_n,\epsilon)$ is a proper holomorphic curve in $B(q_n, \epsilon)$ passing through $q_n$. By Lelong inequality we have area($A_n\cap B(q_n,\epsilon))\geq c\ \epsilon^2$ ($c$ depending only on $X$). This is the contradiction.

\thm{Lelong inequality [\LE]} Let $C$ be a proper holomorphic curve in $B(0,\epsilon)\subset \C^n$, passing through $0$. Then area($C)\geq\pi\ \epsilon^2$. \fthm

\demo{Proof} For $0<r<\epsilon$ put $C_r=C\cap B(0,r)$ and $a(r)=\text{area}(C_r)$. Note that $\partial C_r=C\cap \partial B(0,r)$ by properness. We claim that $\frac {a(r)} {r^2} $ is increasing. Assuming this we get
$$ \frac {\text{area}(C)}{\epsilon^2}=\lim_{r\to \epsilon}\frac {a(r)} {r^2} \geq \lim_{r\to 0}\frac {a(r)} {r^2} \geq \pi$$ 
as $C$ goes through $0$. As for the claim if $0<r<s<\epsilon$ we have $$\frac {a(s)} {s^2}-\frac {a(r)} {r^2}=
\frac 1 {s^2}\int_{C_s}dd^c\Vert z\Vert^2-\frac 1 {r^2}\int_{C_r}dd^c\Vert z\Vert^2=\frac 1 {s^2}\int_{\partial C_s}d^c\Vert z\Vert^2-\frac 1 {r^2}\int_{\partial C_r}d^c\Vert z\Vert^2$$
$$=\int_{\partial C_s}d^c\log\Vert z\Vert^2-\int_{\partial C_r}d^c\log\Vert z\Vert^2=\int_{C_s\setminus C_r}dd^c\log\Vert z\Vert^2\geq0$$
by Stokes theorem and the fact that $dd^c\log\Vert z\Vert^2$ is a positive form.
\enddemo

A consequence is the compactness of holomorphic discs under an area bound in absence of rational curves.
\thm{Compactness of discs} Assume $X$ does not contain any rational curve. Let $f_n :D\to X$ be a sequence of holomorphic discs of bounded area. Then it converges toward a holomorphic disc after extracting. \fthm

\noindent Indeed if $(f_n)$ does not admit any converging subsequence, then $\Vert f_n'\Vert$ has to blow up on a smaller disc (if not apply Ascoli theorem). So by Brody lemma we get an entire curve of bounded area, contradicting the absence of rational curve.

\null
Let us now sketch a proof of our variant under this hypothesis (absence of rational curves) for simplicity.

\null
\noindent {\bf Proof of the variant}. Recall that we are given an Ahlfors current
$T=\lim \frac {[f_n(D)]}{a_n}$ with $\frac {l_n}{a_n}\to 0$, charging $K$. Note that the hypothesis forces $a_n$ to blow up. Let $U$ be a neighbourhood of $K$. For simplicity we will just construct an entire curve passing through $U$ out of these discs. The full statement would follow in the same way taking a sequence of neighbourhoods shrinking to $K$ and extracting diagonally. By assumption for large $n$ we have $$\frac {\text{area}(f_n(D)\cap U )}{a_n}\geq \epsilon.$$
 
\null
The strategy for constructing our entire curve is to look for germs of it, small discs passing through $U$, try and double them as long as we could and pass to the limit using an area bound. Let us make this scheme more precise.

\null
Consider collections $\Cal F_n$ of round discs centered in $D$ depending on $n$. They are {\it disjoint} if the discs in $\Cal F_n$ are disjoint, and {\it consistent} if $\frac {\text{card}(\Cal F_n)}{a_n}$ is bounded from below, say by $\delta$. Let us call disjoint consistent collections simply {\it families}. Families automatically contain subfamilies with an area bound. Indeed for at least half of the discs $d$ of $\Cal F_n$ we have area($f_n(d\cap D))\leq \frac 2 {\delta}$.

\null
Let us explain how to get families of germs through $U$. For each point $z$ of $D$ let $d_z$ be the smallest disc centered at $z$ such that area($f_n(d_z\cap D)\cap U)=1$. By Besicovich lemma [\MA] we may extract from this collection a fixed number $N$ of disjoint subcollections covering $D$. Take the one with most area in $U$. This gives the families $\Cal F_n$ of germs.

\null
We will double them now. Given a disc $d$ write $2d$ for the disc centered at the same point and of radius twice. Consider the collections $2\Cal F_n=\{2d, d \in \Cal F_n\}$. If we are able to extract families from them, then we say that the doubling process works. In this case as already seen we get an area bound.  Moreover we may suppose that the discs $2d$ stay mostly inside $D$. If not we would get families of such discs for which on one hand $2d\cap \partial D$ is relatively big, say $\partial D$ cuts $2d$ roughly by half, and on the other hand length($f_n(2d\cap \partial D)\leq \frac {2l_n}{\delta a_n}$ (same argument as for the area bound). Reparametrizing the discs $d$ by $D$ we would get at the limit a holomorphic map $f :2D^+\to X$ (where $2D^+$ is the upper half of $2D$) such that $f\vert_{[-2,2]}$ is constant (because $\frac {l_n}{a_n}\to 0$) but $f$ is not (because area($f(2D^+)\geq1$). This is impossible.

\null
Assume that we are able to iterate this doubling process indefinitely. Then we have families
$2^k\Cal F_n$ for each $k$. Reparametrizing the $d$ in $\Cal F_n$ by $D$, we deduce that the holomorphic maps $f_n\circ r_n :2^kD\to X$ have area bounds for each $k$. Passing to the limit we get a holomorphic map $f:\C\to X$ such that area($f(\C)\cap \overline U) \geq 1$.

\null
But in general the doubling process may fail at some step, at the first for instance. Then a combinatorial argument shows that there exist in $\Cal F_n$ arbitrary long chains of discs $d_1,...,d_k$ of fast decreasing size, such that $2d_k\subset2d_{k-1}\subset...\subset 2d_2\subset 2d_1$. Now instead of reasoning with discs, we involve the annuli $2d_i\setminus 2d_{i+1}$. Actually we have families of such annuli. With roughly the same arguments as before for the discs, by reparametrizing the middle annulus of these chains and passing to the limit, we get this time a holomorphic map $f:\C^*\to X$. We even have area($f(\C^* )\cap \overline U)=\infty$. See [\DU] for more details.

\null
\noindent {\bf Bloch principle}. Let us mention to finish an application to Bloch principle due to M. McQuillan. In the classical setting Bloch principle reads as follows. Consider a configuration of four lines in general position in $P^2(\C)$. It defines three extra lines joining opposite double points, the {\it diagonals}. Borel [\BO] remarked that the complement of the configuration is hyperbolic {\it modulo the diagonals}. Any entire curve avoiding the configuration has to be contained in a diagonal. Bloch [\BL] translated this in terms of discs. A diverging sequence of holomorphic discs avoiding the configuration has to converge toward the diagonals. In other words diverging sequences of discs localize precisely where entire curves are. In [\MC] McQuillan extends this principle to general logarithmic surfaces by using the variant.

 \Refs

%\medskip
\widestnumber\no{99}
\refno=0
\bref \by L. Ahlfors \paper Zur Theorie der \"Uberlagerungsfl\"achen \jour Acta Math. \vol65\yr1935\pages157--194
\endref
 
\bref \by F. Berteloot, J. Duval \paper Sur l'hyperbolicit\'e de 
certains compl\'ementaires \jour Ens. Math. 
\vol47\yr2001\pages253--267
\endref
 
\bref \by A. Bloch \paper Sur les syst\`emes de fonctions holomorphes \`a vari\'et\'es lin\'eaires lacunaires \jour 
 Ann. Sci. Ecole Norm. Sup. \vol43\yr1926\pages309--362
\endref

\bref \by E. Borel \paper Sur les z\'eros des fonctions enti\`eres \jour Acta Math. \vol20\yr1897\pages357--396
\endref

\bref \by R. Brody \paper Compact manifolds and hyperbolicity \jour 
Trans. Amer. Math. Soc. \vol235\yr1978\pages213--219
\endref

\bref \by C. Ciliberto, M. Zaidenberg \paper Scrolls and hyperbolicity \jour Internat. J. Math. \vol24\yr2013\pages1350026-25p
\endref

\bref \by J. Duval \paper Une sextique hyperbolique dans $P^3(\C)$ \jour Math. Ann. \vol330\yr2004\pages473--476
\endref 

\bref \by J. Duval \paper Sur le lemme de Brody \jour Invent. Math. \vol173\yr2008\pages305--314
\endref
 
 \bref \by M. Green \paper Some Picard theorems for holomorphic maps 
 to algebraic varieties \jour Amer. J.  Math. 
 \vol97\yr1975\pages43--75
 \endref

\bref \by P. Griffiths, J. Harris \book Principles of algebraic geometry \publ Wiley \yr1978\publaddr New York
\endref

\bref \by D. T. Huynh \paper Examples of hyperbolic hypersurfaces of low degree in projective spaces \jour arXiv: 1507.03542, to appear in IMRN
\endref

\bref \by  B. Kleiner \paper Hyperbolicity using minimal surfaces \jour preprint
\endref

\bref \by P. Lelong \paper Propri\'et\'es m\'etriques des vari\'et\'es analytiques complexes d\'efinies par une \'equation\jour Ann. Sci. Ecole Norm. Sup. \vol67 \yr1950 \pages393--419 
\endref  
 
\bref \by P. Mattila \book Geometry of sets and measures in euclidean spaces\publ Cambridge University Press \yr 1995 \publaddr Cambridge
\endref

\bref \by M. McQuillan \paper The Bloch principle \jour arXiv:1209.5402
\endref

\bref \by M. Paun \paper Vector fields on the total space of hypersurfaces in the projective space and hyperbolicity \jour Math. Ann. \vol340\yr2008\pages875--892
\endref 

\bref \by A. Ros \book The Gauss map of minimal surfaces, {\rm in} Differential geometry, Valencia 2001 \pages235--252 \publ World Sci. \yr2002\publaddr River Edge
\endref

\bref \by B. Shiffman, M. Zaidenberg \paper New examples of Kobayashi hyperbolic surfaces in $P^3(\C)$ \jour Funct. Anal. Appl. \vol39\yr2005\pages76--79
\endref

\bref \by J. Winkelmann \paper On Brody and entire curves \jour Bull. Soc. Math. France \vol135\yr2007\pages25--46 
\endref

\endRefs

\enddocument